\input amstex\documentstyle{amsppt}  
\pagewidth{12.5cm}\pageheight{19cm}\magnification\magstep1
\topmatter
\title Twelve bridges from a reductive group to its Langlands dual\endtitle
\author G. Lusztig\endauthor
\address{Department of Mathematics, M.I.T., Cambridge, MA 02139}\endaddress
\thanks{Supported in part by the National Science Foundation}\endthanks
\endtopmatter   
\document

\define\dw{\dot w}

\define\ds{\dot s}

\define\dX{\dot X}

\define\uf{\un f}
\define\ua{\un a}
\define\uv{\un v}
\define\ut{\un t}

\define\uE{\un E}

\define\uU{\un U}

\define\da{\dagger}

\define\Lie{\text{\rm Lie }}

\define\si{\sim}

\define\bK{\bar K}
\define\bC{\bar C}

\define\lb{\linebreak}

\define\bin{\binom}
\define\op{\oplus}

\define\part{\partial}
\define\em{\emptyset}

\define\ra{\rangle}
\define\n{\notin}
\define\iy{\infty}
\define\m{\mapsto}
\define\do{\dots}
\define\la{\langle}

\define\lra{\leftrightarrow}

\define\sm{\smallmatrix}
\define\esm{\endsmallmatrix}
\define\sub{\subset}    

\define\T{\times}
\define\ti{\tilde}
\define\nl{\newline}
\redefine\i{^{-1}}
\define\fra{\frac}
\define\un{\underline}
\define\ov{\overline}
\define\ot{\otimes}
\define\bbq{\bar{\QQ}_l}

\define\Ad{\text{\rm Ad}}
\define\Hom{\text{\rm Hom}}

\define\Aut{\text{\rm Aut}}

\define\Gal{\text{\rm Gal}}

\define\tr{\text{\rm tr}}

\define\a{\alpha}

\redefine\c{\chi}
\define\g{\gamma}
\redefine\d{\delta}
\define\e{\epsilon}
\define\et{\eta}
\define\io{\iota}
\redefine\o{\omega}
\define\p{\pi}
\define\ph{\phi}
\define\ps{\psi}
\define\r{\rho}
\define\s{\sigma}
\redefine\t{\tau}
\define\th{\theta}
\define\k{\kappa}
\redefine\l{\lambda}
\define\z{\zeta}
\define\x{\xi}

\redefine\G{\Gamma}
\redefine\D{\Delta}

\define\Si{\Sigma}

\redefine\L{\Lambda}
\define\Ph{\Phi}
\define\Ps{\Psi}

\define\CC{\bold C}

\define\FF{\bold F}

\define\NN{\bold N}

\define\QQ{\bold Q}
\define\RR{\bold R}

\define\WW{\bold W}
\define\ZZ{\bold Z}

\define\ca{\Cal A}
\define\cb{\Cal B}
\define\cc{\Cal C}
\define\cd{\Cal D}
\define\ce{\Cal E}

\define\ch{\Cal H}
\define\ci{\Cal I}

\define\ck{\Cal K}
\define\cl{\Cal L}
\define\cm{\Cal M}

\define\co{\Cal O}

\define\car{\Cal R}
\define\cs{\Cal S}
\define\ct{\Cal T}

\define\cx{\Cal X}

\define\fa{\frak a}

\define\fm{\frak m}

\define\fr{\frak r}

\define\fA{\frak A}

\define\tg{\ti g}

\define\tm{\ti m}

\define\tI{\ti I}

\define\tW{\ti W}

\define\sh{\sharp}

\define\che{\check}
\define\cha{\che{\a}}

\define\ABV{ABV}
\define\AJS{AJS}
\define\ABG{ABG}
\define\CH{C}
\define\DL{DL}
\define\DR{Dr}
\define\GI{Gi}
\define\HT{HT}
\define\IM{IM}
\define\KW{KW}
\define\KT{KT}
\define\KL{KL1}
\define\KLL{KL2}
\define\KLLL{KL3}
\define\KI{Ki}
\define\KO{Ko}
\define\LF{Lf}
\define\LA{La1}
\define\LAA{La2}
\define\LB{L1}
\define\LC{L2}
\define\LD{L3}
\define\LE{L4}
\define\LG{L5}
\define\LI{L6}
\define\LII{L7}
\define\LJ{L8}
\define\LK{L9}
\define\LM{L10}
\define\LN{L11}
\define\LNN{L12}
\define\LO{L13}
\define\LP{L14}
\define\LQ{L15}
\define\LR{L16}
\define\LS{L17}
\define\LU{L18}
\define\LT{L19}
\define\LV{LV}
\define\MK{MK}
\define\SL{SL}
\define\VO{V1}
\define\VOO{V2}
\define\XI{Xi}
\define\ZE{Ze}

\subhead Introduction\endsubhead
These notes are based on a series of lectures given by the author at the Central China 
Normal University in Wuhan (July 2007). The aim of the lectures was to provide an
introduction to the Langlands philosophy.

According to a theorem of Chevalley, connected reductive split algebraic groups over a 
fixed field $A$ are in bijection with certain combinatorial objects called root systems. 
Now there is a natural involution on the collection of root systems in which roots and 
coroots are interchanged. This corresponds by Chevalley's theorem to an involution on the 
collection of connected reductive split algebraic groups over $A$. The image of a group 
under this involution is called the Langlands dual of that group.

In the 1960's Langlands made the remarkable discovery that some features about the 
representations of a reductive group (such as classification) should be recorded in terms 
of data in the Langlands dual group. He thus formulated two conjectures: one involving
groups over a local field and one involving automorphic representations with respect to a 
group over a global field. 

In these notes we try to give several examples of "bridges" which connect some aspect of
the collection $(G_A)$ of Chevalley groups attached to a root system $\car$ and to various
commutative rings $A$ and some aspect of the analogous collection $(G^*_A)$ of Chevalley 
groups attached to the dual root system $\car^*$ and to various commutative rings $A$. By 
"aspect" we mean something about the structure of one of the groups $G_A$ or of its 
representations or of an associated object such as the affine Hecke algebra $\ch$. In each
case the existence of the bridge is surprising due to the fact that $(G_A)$ and $(G^*_A)$ 
are related only through a very weak connection (via their root systems); in particular 
there is no direct, elementary construction which produces the Langlands dual group from a
given group.

In fact we describe twelve such bridges (some conjectural) the first three of which are
very famous and were found by Langlands himself.

(I) A (conjectural) bridge from irreducible admissible representations of $G_K$ (where $K$
is a finite extension of $\QQ_p$) to certain conjugacy classes of homomorphisms of the Weil
group $W_K$ to $G^*_\CC$. (See \S10.) This bridge contains almost as a special case a
bridge from irreducible representations of $\ch$ specialized at a non-root of $1$ and
conjugacy classes of certain pairs of elements in $G^*_\CC$. (See \S9.)

(II) A bridge from irreducible admissible representations of $G_\RR$ to certain conjugacy 
classes of homomorphisms of the Weil group $W_\RR$ to $G^*_\CC$. (See \S11.)

(III) A (conjectural) bridge connecting certain automorphic representations attached to 
$G_k$ ($k$ a function field over $\FF_p$) and certain homomorphisms of $\Gal(\bar k/k)$ 
into $G^*_\CC$. (See \S12.)

(IV) A bridge from cells in the affine Weyl group constructed from $\ch$ to unipotent 
classes in $G_\CC^*$. (See \S13.)

(V) A bridge from "special unipotent pieces" in $G_{\ov{\FF}_p}$ to "special unipotent 
pieces" in $G^*_{\ov{\FF}_p}$. (See \S14.)

(VI) A bridge from irreducible representations of $G_{\FF_q}$ to certain "special" 
conjugacy classes in $G_\CC^*$. (See \S16.)

(VII) A bridge from character sheaves on $G_\CC$ to "special" conjugacy classes in 
$G_\CC^*$. (See \S17.)

(VIII) A bridge constructed by Vogan connecting certain intersection cohomology spaces
associated to symmetric spaces of $G_\CC$ and similar objects for $G_\CC^*$. (See \S18.)

(IX) A (partly conjectural) bridge connecting multiplicities in standard modules of $G_K$ 
(as in (II)) or $G_\RR$ with intersection cohomology spaces arising from the geometry of
$G^*_\CC$. (See \S19.)

(X) A bridge connecting the tensor product of two irreducible finite dimensional 
representations of $G_\CC$ with the convolution of certain perverse sheaves on the affine
Grassmannian attached to $G^*_{\CC((\e))}$. (See \S20.)

(XI) A bridge connecting the canonical basis of the plus part of the enveloping algebra
attached to $G_\QQ$ with certain subsets of the totally positive part of the upper 
triangular subgroup of $G^*_A$ where $A=\RR[[\e]]$. (See \S21.)

(XII) A (partly conjectural) bridge connecting the characters of irreducible modular 
representations of $G_{\ov{\FF}_p}$ with certain intersection cohomology spaces 
associated with the geometry of $G^*_{\CC((\e))}$. (See \S22.)
\nl
Note that in some cases (such as the very important Case III) our treatment is only a very 
brief sketch. Moreover to simplify the exposition we restrict ourselves to the case of 
split groups.

We now describe the contents of these notes. In \S1 we introduce root systems. In \S2 we
use an idea of McKay (extended by Slodowy) to construct the irreducible simply connected
root systems. In \S3 we introduce the affine Weyl group. In \S4 we introduce the affine 
Hecke algebra and its asymptotic version \cite{\LK}. In \S5 we define the $\ZZ$-form of the
coordinate ring of a Chevalley group. We do not follow the original approach of \cite{\CH}
but rather the approach of Kostant \cite{\KO}. In \S6 we define the Chevalley groups. In 
\S7 we define the Weyl modules. In \S8 we define the Langlands dual group. In \S9-\S22 we 
discuss the various bridges mentioned above.

I wish to thank David Vogan for some useful comments on a first version of these notes.

\subhead 1. Root systems\endsubhead
A {\it root system} is a collection 

$\car=(Y,X,\la,\ra,\cha_i,\a_i(i\in I))$ 
\nl
where $Y,X$ are finitely generated free abelian groups, $\la,\ra:Y\T X@>>>\ZZ$ is a perfect
bilinear pairing, $I$ is a finite set, $\cha_i(i\in I)$ are elements of $Y$ and 
$\a_i(i\in I)$ are elements of $X$ such that $\la\cha_i,\a_i\ra=2$ for all $i$,
$\la\cha_i,\a_j\ra\in-\NN$ for all $i\ne j$; it is assumed that the following (equivalent)
conditions are satisfied:

(i) there exist $c_i\in\ZZ_{>0}(i\in I)$ such that the matrix 
$(c_i\la\cha_i,\a_j\ra)_{i,j\in I}$ is symmetric, positive definite;

(ii) there exist $c'_i\in\ZZ_{>0}(i\in I)$ such that the matrix 
$(\la\cha_i,\a_j\ra c'_j)_{i,j\in I}$ is symmetric, positive definite.
\nl
By the equivalence of (i),(ii), 

$\car^*=(X,Y,\la,\ra',\a_i,\cha_i(i\in I))$ 
\nl
where $\la x,y\ra'=\la y,x\ra$ for $x\in X,y\in Y$, is again a root system, said to be the
dual of $\car$. Note that $\car^{**}=\car$ in an obvious way.

For $\car$ as above let $W$ be the (finite) subgroup of $\Aut(Y)$ generated by the 
automorphisms $s_i:y\m y-\la y,\a_i\ra\cha_i$ ($i\in I$) or equivalently the subgroup of 
$\Aut(X)$ generated by the automorphisms $s_i:x\m x-\la\cha_i,x\ra\a_i$ ($i\in I$); these 
two subgroups may be identified by taking contragredients. We say that $W$ is the {\it Weyl
group} of $\car$; it is also the Weyl group of $\car^*$. Let 

$X^+=\{\l\in X;\la\cha_i,\l\ra\in\NN\text{ for all }i\in I\}$,

$Y^+=\{y\in Y;\la y,\a_i\ra\in\NN\text{ for all }i\in I\}$.
\nl
For $\l,\l'\in X$ write $\l'\ge\l$ if $\l'-\l\in\sum_i\NN\a_i$ and $\l'>\l$ if $\l'\ge\l$,
$\l'\ne\l$.

We say that $\car$ is simply connected if $Y=\sum_i\ZZ\cha_i$. We say that $\car$ is 
adjoint if $X=\sum_i\ZZ\a_i$. We say that $\car$ is semisimple if $X/\sum_i\ZZ\a_i$ is 
finite or equivalently $Y/\sum_i\ZZ\cha_i$ is finite. We say that $\car$ is irreducible if
$I\ne\em$ and there is no partition $I=I'\cup I''$ of $I$ such that $I',I''$ are $\ne\em$
and $\la\cha_i,\a_j\ra=0$ for all $i\in I',j\in I''$. 

\subhead 2. Subgroups of $SL_2(\CC)$ and root systems\endsubhead
In \cite{\MK} McKay discovered a remarkable direct connection between finite subgroups of
$SL_2(\CC)$ and "simply laced affine Dynkin diagrams". Slodowy \cite{\SL} extended this to
a connection between certain pairs of subgroups of $SL_2(\CC)$ (one contained in the other)
and "affine Dynkin diagrams". 

Let $\G,\G'$ be two finite subgroups of $SL_2(\CC)$ such that $\G$ is a normal subgroup of
$\G'$ with $\G'/\G$ cyclic. We show how to attach to $(\G,\G')$ a root system (we use an
argument generalizing one in \cite{\LNN, 1.2}). Let $\cx$ be the category of finite 
dimensional complex representations of $\G$ which can be extended to 
representations of $\G'$. Let $\r_i(i\in\tI)$ be the indecomposable objects of $\cx$ up to
isomorphism that is, the representations of $\G$ which are restrictions of irreducible 
representations of $\G'$. Let $i_0\in\tI$ be such that $\r_{i_0}$ is the trivial 
representation of $\G$ on $\CC$. Let $\s$ be the obvious representation of $\G$ on $\CC^2$;
we have $\s\in\cx$. Let $V$ be the $\RR$-vector space with basis $\{i;i\in\tI\}$. Any 
object $\r$ of $\cx$ gives rise to a vector $\un\r=\sum_in_ii\in V$ where 
$\r\cong\op_i\r_i^{\op n_i}$; here $n_i\in\NN$. For $i\in\tI$ we have $\r_i\ot\s\in\cx$ and
$\un{\r_i\ot\s}=\sum_{j\in\tI}c_{ij}j$ with $c_{ij}\in\NN$. 

For $i\in\tI$, $\r_i$ is the direct sum of $m_i$ irreducible representations of $\G$ (each 
with multiplicity $1$). Let $[,]$ be the bilinear form on $V$ with values in $\RR$ given by
$[i,j]=(2\d_{ij}-c_{ij})m_j$ for $i,j\in\tI$. Let $x=\sum_ix_ii\in V$, 
$x'=\sum_ix'_ii\in V$ where $x_i,x'_i\in\RR$. For $\g\in\G$ let $\l_g$ be an eigenvalue of
$\g$ on $\CC^2$. We have
$$\align&[x,x']=|\G|\i\sum_{i,j;\g\in\G}x_ix'_j
\tr(\g,\r_i)\ov{\tr(\g,\r_j)}(2-\l_g-\ov{\l_g})\\&=|\G|\i\sum_{\g\in\G}
(\sum_ix_i\tr(\g,\r_i)|1-\l_\g|)(\sum_jx'_j\ov{\tr(\g,\r_j)}|1-\l_\g|).\endalign$$
In particular,
$$[x,x]=|\G|\i\sum_{\g\in\G}|\sum_ix_i\tr(\g,\r_i)|1-\l_\g||^2\ge0.$$
If $[x,x]=0$ then for any $\g\in\G$ we have $\sum_ix_i\tr(\g,\r_i)|1-\l_\g|=0$ that is, for
any $\g\in\G-\{1\}$ we have $\sum_ix_i\tr(\g,\r_i)=0$ that is, there exists $c\in\RR$ such
that $x=c\un{\fr}$ where $\fr\in\cx$ is the regular representation of $\G$; if in addition
we have $x_{i_0}=0$ then we see that $c=0$ hence $x=0$.

Let $I=\tI-\{i_0\}$. Let $Y$ be the subgroup of $V$ generated by $\{i;i\in I\}$. For 
$i\in I$ let $\cha_i=i\in Y$. Let $X=\Hom(Y,\ZZ)$. Let $\la,\ra:Y\T X@>>>\ZZ$ be the 
obvious pairing. For $j\in I$ define $\a_j\in X$ by $\la\cha_i,\a_j\ra=(i,j)$. We have 
$\la\cha_i,\a_i\ra=2$ for all $i\in I$ and $\la\cha_i,\a_j\ra=-c_{ij}\in-\NN$ for $i\ne j$
in $I$. By the argument above the matrix $(\la\cha_i,\a_j\ra m_j)_{i,j\in I}$ is symmetric
and positive definite. Hence $\car=(Y,X,\la,\ra,\cha_i,\a_i(i\in I))$ is a (simply 
connected) root system.

Note that all simply connected irreducible root systems are obtained by this construction
exactly once (up to isomorphism) from pairs $(\G,\G')$ as above (up to conjugacy) with 
$\G\ne\{1\}$ and with the property that any element of $\G'$ which commutes with any 
element of $\G$ is contained in $\G$. Such pairs are classified as follows:

(a) $\G=\G'$ is a cyclic group $\ZZ_n$ of order $n\ge2$;

(b) $\G=\G'$ is a binary dihedral group $\cd_{4n}$ of order $4n\ge8$;

(c) $\G=\G'$ is a binary tetrahedral group $G_{24}$ of order $24$;

(d) $\G=\G'$ is a binary octahedral group $G_{48}$ of order $48$;

(e) $\G=\G'$ is a binary icosahedral group $G_{120}$ of order $120$;

(f) $\G'=\cd_{4n},\G=\ZZ_{2n}$ with $n\ge2$;

(g) $\G'=\cd_{8n}$, $\G=\cd_{4n}$ with $n\ge2$;

(h) $\G'=G_{48}$, $\G=G_{24}$;

(i) $\G'=G_{24},\G=\cd_8$.

\subhead 3. Affine Weyl group\endsubhead
Let $\car=(Y,X,\la,\ra,\cha_i,\a_i(i\in I))$ be a root system. Let $W$, $s_i(i\in I)$ be 
as in \S1. 
Let $\tW$ be the semidirect product $W\cdot Y$. We have $\tW=\{wa^y;w\in W,y\in Y\}$ where
$a$ is a symbol; the multiplication is given by $(wa^y)(w'a^{y'})=ww'a^{w\i(y)+y'}$ for 
$w,w'\in W$, $y,y'\in Y$. We identify $Y$ with its image under the homomorphism $y\m1a^y$ 
(a normal subgroup of $\tW$) and $W$ with its image under the homomorphism $w\m wa^0$.

Let $R$ be the set of elements of $X$ of the form $w(\a_i)$ for some $w\in W,i\in I$. Let 
$\che R$ be the set of elements of $Y$ of the form $w(\cha_i)$ for some $w\in W,i\in I$.
There is a unique $W$-equivariant bijection $\a\lra\cha$ between $R$ and $\che R$ such that
$\a_i\lra\cha_i$ for any $i\in I$. For $\a\in R$ we set $s_\a=ws_iw\i\in W$ where 
$\a=w(\a_i),w\in W,i\in I$. Note that $s_\a$ is well defined. Let $R_{min}$ be the set of 
all $\a\in R$ such that the following holds: if $\a'\in R,\a'\le\a$ then $\a'=\a$.

Let $R^+=R\cap(\sum_i\NN\a_i)$, $R^-=-R^+$. We have $R=R^+\cup R^-$.
Following Iwahori and Matsumoto \cite{\IM}, we define a function $l:\tW@>>>\NN$ by
$$l(wa^y)
=\sum_{\a\in R^+;w(\a)\in R^-}|\la y,\a\ra+1|+\sum_{\a\in R^+;w(\a)\in R^+}|\la y,\a\ra|.$$
Let $S$ be the subset of $W$ consisting of the involutions $s_i(i\in I)$ and the 
involutions $s_\a a^{\cha}$ with $\a\in R_{min}$. Note that $l|_S=1$.

If $y\in Y^+$ we have $l(a^y)=\sum_{\a\in R^+}\la y,\a\ra$. Hence for $y,y'\in Y^+$ we have
$l(a^y\cdot a^{y'})=l(a^{y+y'})=l(a^y)+l(a^{y'})$.

\subhead 4. Affine Hecke algebra\endsubhead
We preserve the notation of \S3. Let $\ca=\ZZ[v,v\i]$ where $v$ is an indeterminate. Let 
$\ch$ be the associative $\ca$-algebra with $1$ with generators $T_w (w\in\tW)$ and 
relations 

$(T_s-v)(T_s+v\i)=0$ for $s\in S$,

$T_wT_{w'}=T_{ww'}$ for $w,w'\in\tW$ such that $l(w)+l(w')=l(ww')$.
\nl
We have $T_1=1$ and $\{T_w;w\in\tW\}$ is an $\ca$-basis of $\ch$. Let $h\m h^\da$ be the 
$\ca$-algebra involution of $\ch$ such that $T_w^\da=(-1)^{l(w)}T_w\i$. Let $h\m\bar h$ be
the ring involution of $\ch$ such that $\ov{T_w}=T_{w\i}\i$ for $w\in\tW$ and 
$\ov{v^n}=v^{-n}$ for $n\in\ZZ$. Let $z\in\tW$. According to \cite{\KL} there is a unique 
element $c_z\in\ch$ such that $\ov{c_z}=c_z$ and $c_z=\sum_{w\in\tW}p_{w,z}T_w$ where 
$p_{w,z}\in v\i\ZZ[v\i]$ for all $w\ne z$, $p_{z,z}=1$ and $p_{w,z}=0$ for all but finitely
many $w$. Note that $\{c_w;w\in\tW\}$ is an $\ca$-basis of $\ch$. For $x,y\in\tW$ we write
$c_xc_y=\sum_{z\in\tW}h_{x,y,z}c_z$ where $h_{x,y,z}\in\ca$ is $0$ for all but finitely 
many $z$. Let $z\in\tW$. According to \cite{\LK} there is a unique $a(z)\in\NN$ such that 
$h_{x,y,z}\in v^{a(z)}\ZZ[v\i]$ for any $x,y$ and $h_{x,y,z}\n v^{a(z)-1}\ZZ[v\i]$ for some
$x,y$. We have $h_{x,y,z}=\g_{x,y,z\i}v^{a(z)}\mod v^{a(z)-1}\ZZ[v\i]$ where 
$\g_{x,y,z\i}\in\ZZ$. Let $J$ be the free abelian group with basis $\{t_w;w\in\WW\}$. 
Consider the $\ZZ$-algebra structure on $J$ such that 
$t_xt_y=\sum_{z\in\tW}\g_{x,y,z\i}t_z$ for all $x,y$ in $\tW$. This structure is 
associative and the ring $J$ has a unit element $1$ of the form $\sum_{d\in\cd}t_d$ where 
$\cd$ is a finite subset of $\tW$ consisting of involutions \cite{\LK}. For $x,y$ in $\tW$
we write $x\si y$ when $t_xt_ut_y\ne0$ for some $u\in\tW$. This is an 
equivalence relation on $\tW$; the equivalence classes are called two-sided cells. For any
two-sided cell $c$ let $J_c$ be the subgroup of $J$ spanned by $\{t_z;z\in c\}$. This is a
subring of $J$ with unit $\sum_{d\in\cd\cap c}t_d$. We have $J=\op_cJ_c$ as rings. Consider
the $\ca$-linear map $\ph:\ch@>>>\ca\ot J$ given by 
$$\ph(c_x^\da)=\sum_{z\in\tW,d\in\cd;a(d)=a(z)}h_{x,d,z}t_z$$
for any $x\in\tW$. This is an $\ca$-algebra homomorphism \cite{\LK}.

\subhead 5. Coordinate ring\endsubhead
Let $\car=(Y,X,\la,\ra,\cha_i,\a_i(i\in I))$ be a root system. Let $\uf$ be the associative
$\QQ$-algebra with $1$ defined by the generators $\th_i(i\in I)$ and the Serre relations
$$\sum_{a,b\in\NN;a+b=1-\la\cha_i,\a_j\ra}(-1)^a(\th_i^a/a!)\th_j(\th_i^b/b!)$$
for $i\ne j$ in $I$. Let ${}^0\uU$ be the symmetric algebra of $\QQ\ot Y$. Let $\uU$ be the
$\QQ$-algebra with $1$ defined by the generators $x^+,x^-(x\in\uf)$, $\ua\in{}^0\uU$ and 
the relations:

$x\m x^+$ is an algebra homomorphism $\uf@>>>\uU$ respecting $1$;

$x\m x^-$ is an algebra homomorphism $\uf@>>>\uU$ respecting $1$;

$\ua\m\ua$ is an algebra homomorphism ${}^0\uU@>>>\uU$ respecting $1$;

$y\th_i^+-\th_i^+y=\la y,\a_i\ra\th_i^+$ for $y\in Y,i\in I$;

$y\th_i^--\th_i^-y=-\la y,\a_i\ra\th_i^-$ for $y\in Y,i\in I$;

$\th_i^+\th_j^--\th_j^-\th_i^+=\d_{ij}\cha_i$ for $i,j$ in $I$.
\nl
Define a $\QQ$-algebra homomorphism $\D:\uU@>>>\uU\ot\uU$ by
$\D(\th_i^+)=\th_i^+\ot1+1\ot\th_i^+$, $\D(\th_i^-)=\th_i^-\ot1+1\ot\th_i^-$ for $i\in I$,
$\D(y)=y\ot1+1\ot y$ for $y\in Y$. Define a $\QQ$-algebra isomorphism 
$\Si:\uU@>>>\uU^{opp}$
by $\Si(\th_i^+)=-\th_i^+$, $\Si(\th_i^-)=-\th_i^-$ for $i\in I$, $\Si(y)=-y$ for $y\in Y$.

Let $\uf_\ZZ$ be the subring of $\uf$ generated by the elements 
$\th_i^{(n)}:=\th_i^n/n!$, ($i\in I,n\in\NN$). Note that $\uf_\ZZ$ is a lattice in the 
$\QQ$-vector space $\uf$. Following \cite{\KO} we define $\uU_\ZZ$ to be the subring of
$\uU$ generated by the elements $x^+,x^-$ ($x\in\uf_\ZZ$) and \lb
$\bin{y}{n}:=\fra{y(y-1)\do(y-n+1)}{n!}$, ($y\in Y,n\in\NN$). Note that $\uU_\ZZ$ is a 
lattice in the $\QQ$-vector space $\uU$. Hence $\uU_\ZZ\ot_\ZZ\uU_\ZZ$ is a lattice in 
$\uU\ot_\QQ\uU$. Note that $\D$ restricts to a ring homomorphism 
$\uU_\ZZ@>>>\uU_\ZZ\ot_\ZZ\uU_\ZZ$ denoted again by $\D$. Also $\Si$ restricts to a ring 
isomorphism $\uU_\ZZ@>>>\uU_\ZZ^{opp}$ denoted again by $\Si$. 

For any $\uU$-module $M$ and any $x\in X$ we set

$M^x=\{m\in M;ym=\la y,x\ra m\text{ for any }y\in Y\}$.
\nl
Let $\cc$ be the category whose objects are $\uU$-modules $M$ with $\dim_\QQ M<\iy$ such
that $M=\op_{x\in X}M^x$. For any $\QQ$-vector space $V$ we set $V^\da=\Hom_\QQ(V,\QQ)$.
For $M\in\cc$ we define $c_M:M\ot M^\da@>>>\uU^\da$ by $m\ot\x\m[u\m\x(um)]$. Let 
$$\co=\sum_{M\in\cc}c_M(M\ot M^\da),$$
a $\QQ$-subspace of $\uU^\da$. (This agrees with the definition in \cite{\KO} when $\car$
is semisimple.) For $f,f'\in\co$ we define $ff':\uU@>>>\QQ$ by $u\m\sum_sf(u_s)f'(u'_s)$ 
where $\D(u)=\sum_su_s\ot u'_s$, $u_s,u'_s\in\uU$. We have $ff'\in\co$. This defines a 
structure of associative, commutative algebra on $\co$. This algebra has a unit element: 
the algebra homomorphism $\uU@>>>\QQ$ such that $\th_i^+\m0$, $\th_i^-\m0$ for $i\in I$, 
$y\m0$ for $y\in Y$ and $1\m1$. For $f\in\co$ we define a linear function 
$\d(f):\uU\ot\uU@>>>\QQ$ by $u_1\ot u_2\m f(u_1u_2)$. Note that $\co\ot\co$ is naturally a
subspace of $(\uU\ot\uU)^\da$ and that the image of $\d:\co@>>>(\uU\ot\uU)^\da$ is 
contained in the subspace $\co\ot\co$ so that $\d$ defines a linear map $\co@>>>\co\ot\co$
denoted again by $\d$. This is an algebra homomorphism. For $f\in\co$ we define a linear 
function $\s(f):\uU@>>>\QQ$ by $u\m f(\Si(u))$. We have $\s(f)\in\co$; thus $\s$ defines a
linear map $\co@>>>\co$ denoted again by $\s$. This is an algebra homomorphism. Define 
$\e:\co@>>>\QQ$ by $f\m f(1)$. Note that the commutative algebra $\co$ with the 
comultiplication $\d$, the antipode $\s$ and the counit $\e$ is a Hopf algebra over $\QQ$.

Let $f\in\co$. There is a unique collection $(f^x)_{x\in X}$ of numbers in $\QQ$ such that 
$f^x=0$ for all but finitely many $x\in X$ and 
$$f(y_1^{n_1}y_2^{n_2}\do y_r^{n_r})=\sum_{x\in X}\la y_1,x\ra^{n_1}\la y_2,x\ra^{n_2}\do 
\la y_r,x\ra^{n_r}f^x$$
for any $y_1,y_2,\do,y_r$ in $Y$ and $n_1,n_2,\do,n_r$ in $\NN$. For example if 
$f=c_M(m\ot\x)$ where $M\in\cc,m\in M,\x\in M^\da$, we have $f^x=\x(m_x)$ where 
$m_x\in M^x$ are defined by $m=\sum_{x\in X}m_x$. 
Note that for any $x\in X$, $f\m f^x$ is a linear function $\co@>>>\QQ$.

Let $\co_\ZZ=\{f\in\co;f(\uU_\ZZ)\sub\ZZ\}$. (This agrees with the definition in \cite{\KO}
when $\car$ is semisimple.) Note that $\co_\ZZ$ is a subring of $\co$. One can show that
$\co_\ZZ$ is a lattice in the $\QQ$-vector space $\co$. Hence $\co_\ZZ\ot_\ZZ\co_\ZZ$ is a
lattice in the $\QQ$-vector space $\co\ot\co$. Note that $\d:\co@>>>\co\ot\co$ restricts to
a ring homomorphism $\d_\ZZ:\co_\ZZ@>>>\co_\ZZ\ot_\ZZ\co_\ZZ$; $\s:\co@>>>\co$ restricts to
a ring isomorphism $\s_\ZZ:\co_\ZZ@>>>\co_\ZZ$; $\e:\co@>>>\QQ$ restricts to a ring 
homomorphism $\e_\ZZ:\co_\ZZ@>>>\ZZ$. The commutative ring $\co_\ZZ$ together with the 
comultiplication $\d_\ZZ$, the antipode $\s_\ZZ$ and the counit $\e_\ZZ$ is a Hopf ring. 
For any $x\in X$ and $f\in\co_\ZZ$ we have $f^x\in\ZZ$.

For any commutative ring $A$ with $1$ we set $\co_A=A\ot\co_\ZZ$. By extension of scalars,
from $\d_\ZZ,\s_\ZZ,\e_\ZZ$ we get $A$-algebra homomorphisms 
$\d_A:\co_A@>>>\co_A\ot_A\co_A$, $\s_A:\co_A@>>>\co_A$, $\e_A:\co_A@>>>A$. The commutative 
$A$-algebra $\co_A$ together with the comultiplication $\d_A$, the antipode $\s_A$ and the
counit $\e_A$ is a Hopf algebra over $A$. For any $x\in X$ the homomorphism 
$\co_\ZZ@>>>\ZZ$ gives rise by extension of scalars to an $A$-linear map $\co_A@>>>A$
denoted again by $f^x$. 

The following two properties are proved in \cite{\LT}:

(i) the $A$-algebra $\co_A$ is finitely generated;

(ii) if $A$ is an integral domain then $\co_A$ is an integral domain.

\subhead 6. Chevalley groups\endsubhead
We preserve the notation of \S5. Let $W,s_i$ be as in \S1. Let $A$ be a commutative ring 
with $1$. As in \cite{\KO} we define $G_A$ to be the 
set of $A$-algebra homomorphisms $\co_A@>>>A$ respecting $1$. For $g,g'\in G_A$ we define
$gg':\co_A@>>>A$ by $f\m\sum_sg(f_s)g'(f'_s)$ where $\d_A(f)=\sum_sf_s\ot f'_s$ with
$f_s,f'_s$ in $\co_A$. Then $gg'\in G_A$ and $(g,g')\m gg'$ is a group structure on $G_A$ 
with unit element $\e_A$. We say that $G_A$ is the {\it Chevalley group} attached to the 
root system $\car$ and to the commutative ring $A$. We write also $G_A^{\car}$ instead of
$G_A$ when we want to emphasize the dependence on $\car$.

If $\k:A@>>>A'$ is a homomorphism of commutative rings with $1$ and $g:\co_A@>>>A$ is in 
$G_A$ then applying to $g$ the functor $A'\ot_A?$ (where $A'$ is regarded as an $A$-algebra
via $\k$) we obtain an $A'$-algebra homomorphism $\co_{A'}@>>>A'$ respecting $1$ which is 
denoted by $\tg$. Now $g\m\tg$ is a group homomorphism $G_A@>>>G_{A'}$ said to be induced 
by $\k$.

For any $i\in I,b\in A$ we define an $A$-linear map $x_i(b):\co_A@>>>A$ by
$$\sum_sa_s\ot f_s\m\sum_s\sum_{n\in\NN}a_sb^nf_s((\th_i^{(n)})^+).$$
Here $a_s\in A,f_s\in\co_\ZZ$. Since $f_s\in\co$ we have $f_s((\th_i^{(n)})^+)=0$ for large
enough $n$ so that the last sum is finite. From the definitions we see that $x_i(b)\in G_A$
and that $b\m x_i(b)$ is an (injective) group homomorphism $A@>>>G_A$.

Similarly, for any $i\in I,b\in A$ we define an $A$-linear map $y_i(b):\co_A@>>>A$ by
$$\sum_sa_s\ot f_s\m\sum_s\sum_{n\in\NN}a_sb^nf_s((\th_i^{(n)})^-).$$
Here $a_s\in A,f_s\in\co_\ZZ$. Again the last sum is finite. From the definitions we see 
that $y_i(b)\in G_A$ and that $b\m y_i(b)$ is an (injective) group homomorphism $A@>>>G_A$.

Let $A^*$ be the group of units of $A$. Let $t=\sum_ra_r\ot y_r\in A^*\ot Y$ with
$a_r\in A^*,y_r\in Y$. We define $\ut:\co_A@>>>A$ by 
$$f\m\sum_r\sum_{x\in X}a_r^{\la y_r,x\ra}f^x.$$
(Since $a_r\in A^*$, $a_r^n$ is defined for any $n\in\ZZ$.) From the definitions we see 
that $\ut\in G_A$ and that $t\m\ut$ is an injective group homomorphism $A^*\ot Y@>>>G_A$
with image denoted by $T_A$ (a commutative subgroup of $G_A$). We identify $A^*\ot Y$ with
$T_A$ via this homomorphism. We write also $T_A^{\car}$ instead of $T_A$ when we want to 
emphasize the dependence on $\car$.

For $i\in I$ we define an element $\ds_i\in G_A$ by $\ds_i=x_i(1)y_i(-1)x_i(1)$. We have 
$\ds_i^2=(-1)\ot\cha_i\in A^*\ot Y=T_A$. Moreover, for $i\ne j$ we have 
$\ds_i\ds_j\ds_i\do=\ds_j\ds_i\ds_j\do$ (both sides have $n$ factors where $n$ is the 
order of $s_is_j$ in $W$). 
It follows that for any $w\in W$ there is a well defined element $\dw\in G_A$ 
such that $\dw=\ds_{i_1}\ds_{i_2}\do\ds_{i_r}$ whenever $w=s_{i_1}s_{i_2}\do s_{i_r}$ with
$r=l(w)$. Note that $\dw T_A\dw\i=T_A$. More precisely, for $t\in T_A$ we have 
$\dw t\dw\i=w(t)$ where $w:t\m w(t)$ is the $W$-action on $T_A$ given by 
$a\ot y\m a\ot w(y)$ for $a\in A^*$, $y\in Y$.

For any sequence $i_1,i_2,\do,i_r$ in $I$ such that $l(s_{i_1}s_{i_2}\do s_{i_r})=r=|R^+|$,
the map $A^r@>>>G_A$ given by
$$\align&(a_1,a_2,\do,a_r)\m \\&x_{i_1}(a_1)\ds_{i_1}x_{i_2}(a_2)\ds_{i_1}\i\do 
\ds_{i_1}\ds_{i_2}\do\ds_{i_{r-1}}x_{i_r}(a_r)\ds_{i_{r-1}}\i\do\ds_{i_2}\i\ds_{i_1}\i
\endalign$$
is injective and its image is a subgroup $U_A^+$ of $G_A$ independent of the choice of
$i_1,i_2,\do,i_r$. (See \cite{\LT}.)

Similarly, for any sequence $i_1,i_2,\do,i_r$ in $I$ such that
$l(s_{i_1}s_{i_2}\do s_{i_r})=r=|R^+|$, the map $A^r@>>>G_A$ given by
$$\align&(a_1,a_2,\do,a_r)\m \\&y_{i_1}(a_1)\ds_{i_1}y_{i_2}(a_2)\ds_{i_1}\i\do 
\ds_{i_1}\ds_{i_2}\do\ds_{i_{r-1}}y_{i_r}(a_r)\ds_{i_{r-1}}\i\do\ds_{i_2}\i\ds_{i_1}\i
\endalign$$
is injective and its image is a subgroup $U_A^-$ of $G_A$ independent of the choice of
$i_1,i_2,\do,i_r$. 

The subgroups $U_A^+,U_A^-$ are normalized by $T_A$. We set $B_A^+=U_A^+T_A=T_AU_A^+$,
$B_A^-=U_A^-T_A=T_AU_A^-$. 

If $A$ is a field, we have a partition $G_A=\cup_{w\in W}B_A^+\dw B_A^+$.

\subhead 7. Weyl modules\endsubhead
We preserve the notation of \S6. 
For any $\l\in X^+$ let $\ct_\l=\sum_i\uf\th_i^{\la\cha_i,\l\ra+1}$, a left ideal of $\uf$;
let $\L_\l=\uf/\ct_\l$, a finite dimensional $\QQ$-vector space. Let $\et$ be the image of
$1\in\uf$ in $\L_\l$. We regard $\L_\l$ as a $\uU$-module in which $x^-$ acts as left
multiplication by $x (x\in\uf)$; $\th_i^+\et=0$ for $i\in I$; $y\et=\la y,\l\ra\et$ for 
$y\in Y$. We say that $\L_\l$ is a {\it Weyl module}. We have $\L_\l\in\cc$. 

For $\l\in X^+$ let 
$$\ct_{\l,\ZZ}=\sum_{i,n;n\ge\la\cha_i,\l\ra+1}\uf\th_i^{(n)}=\ct_\l\cap\uf_\ZZ,$$
a left ideal of $\uf_\ZZ$. Let $\L_{\l,\ZZ}=\uf_\ZZ/\ct_{\l,\ZZ}$. Then $\L_{\l,\ZZ}$ is a
lattice in the $\QQ$-vector space $\L_\l$ and a $\uU_\ZZ$-submodule of $\L_\l$. For a 
commutative ring $A$ with $1$ we set $\L_{\l,A}=\L_{\l,\ZZ}\ot A$.

We write $\co^{opp},\co^{opp}_A$ for $\co,\co_A$ with the opposite comultiplication. Define
$\Xi:\L_\l@>>>\co\ot\L_\l$ by 
$e\m\sum_jc_{\L_l}(e\ot\x'_j)\ot e_j$ where $(e_j)$ is a $\QQ$-basis of $\L_\l$ and 
$(\x_j)$ is the dual basis of $\L_\l^\da$. This makes $\L_\l$ into a $\co^{opp}$-comodule.
Now $\co_\ZZ\ot_\ZZ\L_{\l,\ZZ}$ is a lattice in $\co\ot\L_\l$ and $\Xi$ restricts to 
$\Xi_\ZZ:\L_{\l,\ZZ}@>>>\co_\ZZ\ot_\ZZ\L_{\l,\ZZ}$. By extension of scalars we obtain an 
$A$-linear map $\Xi_A:\L_{\l,A}@>>>\co_A\ot_A\L_{\l,A}$ making $\L_{\l,A}$ into a 
$\co_A^{opp}$-comodule. For any $g:\co_A@>>>A$ which is in $G_A$, we define 
$\r_g:\L_{\l,A}@>>>\L_{\l,A}$ by $e\m\sum_hg(f_h)\ot e_h$ where 
$\Xi_A(e)=\sum_h f_h\ot e_h$, $f_h\in\co_A,e_h\in\L_{\l,A}$. Note that $g\m\r_g$ is a group
action. Thus $\L_{\l,A}$ is a $G_A$-module.

\subhead 8. The Langlands dual group\endsubhead
If $A$ is an algebraically closed field then $G_A=G_A^{\car}$ is a reductive connected 
algebraic group over $A$ with coordinate ring $\co_A$. Moreover, according to Chevalley,
$\car\m G_A^{\car}$ is a bijection 
$$\align&\{\text{root systems up to isom.}\}@>\si>>\\&
\{\text{reductive connected algebraic groups over $A$ up to isom.}\}.\endalign$$

An element of $G_A$ is said to be semisimple if it is conjugate to an element in $T_A$. An
element of $G_A$ is said to be unipotent if it is conjugate to an element in $U_A^+$.

(a) {\it In the remainder of these notes (except in \S15) we fix a root system}
$\car=(Y,X,\la,\ra,\cha_i,\a_i(i\in I))$. Define $W,\tW,\ch$ in terms of $\car$ as in 
\S1,\S3,\S4. Define $\tW^*,\ch^*$ like $\tW,\tW^*$ but in terms of $\car^*$ instead of 
$\car$.

For a commutative ring $A$ with $1$ we write $G_A,T_A$ instead of $G_A^{\car},T_A^{\car}$
and 
$G^*_A,T^*_A$ instead of $G_A^{\car^*},T_A^{\car^*}$. 
We say that $G^*_A$ is the Langlands dual group to $G_A$.
Let $\cb_A$ be the set of subgroups 
of $G_A$ that are conjugate to $B_A^+$ (or equivalently to $B_A^-$). Let $\cb^*_A$ be the 
analogous set defined in terms of $\car^*$ instead of $\car$. For $g\in G^*_A$ we set 
$\cb^*_{A,g}=\{B\in\cb^*_A;g\in B\}$. 

\subhead 9. Representations of affine Hecke algebras\endsubhead
Assume that $X/\sum_i\ZZ\a_i$ has no torsion. We fix $\uv\in\CC^*$. Let 
$\ch_{\uv}=\CC\ot_\ca\ch$ where $\CC$ is regarded as an $\ca$-algebra via $v\m\uv$. Let 
$\Ps_{\uv}$ be the set of all pairs $(s,u)$ up to $G_\CC^*$-conjugation) where 
$s\in G_\CC^*$ is semisimple, $u\in G_\CC^*$ is unipotent and $sus\i=u^{\uv^2}$. 

Assume that $\uv$ is not a root of $1$. The Deligne-Langlands conjecture states that there
is a canonical finite to one surjective map
$$\{\text{irred. $\ch_{\uv}$-modules up to isom.}\}@>>>\Ps_{\uv}.\tag a$$
A refinement of this conjecture was stated in \cite{\LG} namely that the fibre of (a) at
$(s,u)$ should be in natural bijection with the set of irreducible representations (up to
isomorphism) of the group of connected components of the centralizer of $(s,u)$ in 
$G_\CC^*$ which appear in the natural representation of this group on the cohomology of
$\{B\in\cb^*_\CC;s\in B,u\in B\}$. This came from a study of examples connected with
"subregular" unipotent elements in $G_\CC^*$. In \cite{\LII} it was shown that $\ch_{\uv}$
acts naturally on the equivariant $K$-theory of $\cb^*_\CC$ where the parameter of the 
Hecke algebra comes from equivariance with respect to a $\CC^*$-action. In \cite{\LII} it 
was also suggested that one should construct representations of $\ch_{\uv}$ using the 
equivariant $K$-theory of the varieties $\cb^*_{\CC,u}$ for $u\in G_\CC^*$ unipotent. This
was established in \cite{\KLL} which gave a proof of (a) (in the refined form).

In \cite{\XI} it is shown (by a reduction to \cite{\KLL}) that a statement analogous to (a)
(in the refined form) holds also when $\uv$ is allowed to be a root of $1$ in the 
complement of a specific finite set of roots of $1$ depending on $\car$. 

\subhead 10. $p$-adic groups\endsubhead
Let $K$ be a finite extension of the field of $p$-adic numbers ($p$ a prime number). Let 
$\fa$ be the integral closure in $K$ of the ring of $p$-adic integers and let $\fm$ be the
unique maximal ideal of $\fa$ so that $\fa/\fm$ is a finite field with $q$ elements. Let 
$\bK$ be an algebraic closure of $K$. Let $\bar\fa$ be the integral closure of $\fa$ in
$\bK$ and let $\bar\fm$ be the unique maximal ideal of $\bar\fa$ so that $\bar\fa/\bar\fm$
is an algebraic closure of $\fa/\fm$. 
Let $W_K$ be the Weil group of $K$ that is the inverse image under the natural homomorphism
$\p:\Gal(\bK/K)@>>>\Gal(\bar\fa/\bar\fm,\fa/\fm)$ of the subgroup $\ZZ$ of 
$\Gal(\bar\fa/\bar\fm,\fa/\fm)$ consisting of the integer powers of the automorphism 
$x\m x^q$. Let $\ci$ be the kernel of $\p$. We have an exact sequence 
$1@>>>\ci@>>>W_K@>\o>>\ZZ@>>>1$ where $\o$ is the restriction of $\p$. A homomorphism 
$\r:W_K@>>>G_\CC^*$ is said to be admissible if $\r(\ci)$ is finite and $\r(\ph)$ is 
semisimple in $G_\CC^*$ for some/any $\ph\in\o\i(1)$. Let $\Ph^K(G_\CC^*)$ be the set of 
all pairs $(\r,u)$ (up to $G^*_\CC$-conjugacy) where $\r:W_K@>>>G_\CC^*$ is an admissible 
homomorphism and $u\in G_\CC^*$ is a unipotent element such that 
$\r(w)u\r(w)\i=u^{q^{\o(w)}}$ for any $w\in W_K$. We regard $G_K$ as a topological group 
with the $p$-adic topology. An 
irreducible representation $G_K@>>>GL(E)$ (where $E$ is a $\CC$-vector space) is said to be
admissible if the stabilizer of any vector of $E$ is open in $G_K$ and if for any open 
subgroup $H$ of $G_K$ the space of $H$-invariant vectors in $E$ has finite dimension. 
According to the local Langlands conjecture there is a canonical finite to one surjective 
map
$$\{\text{irred. admissible representations of $G_K$ up to isom.}\}@>>>\Ph^K(G_\CC^*).
\tag a$$
This is known to be true in the case where $G_K=GL_n(K)$, see \cite{\HT}. In the general
case but assuming that $\car$ is adjoint, a class of irreducible admissible representations
(called "unipotent") has been described in \cite{\LQ} where a canonical finite to one 
surjective map
$$\{\text{unipotent representations of $G_K$ up to isom.}\}@>>>\Ph^K_1(G_\CC^*)\tag b$$
was constructed; here $\Ph^K_1(G_\CC^*)=\{(\r,u)\in\Ph(G_\CC^*);\r(\ci)=\{1\}\}$. Note that
$\Ph^K_1(G_\CC^*)$ may be identified with the set $\Ps_{\sqrt q}$ in \S9 
and that (b) constitutes a verification of (a) in a special case. Note that some of the
unipotent representations can be understood by the method described in \S9; to understand
the remaining ones one needs the theory of character sheaves and a geometric construction
of certain affine Hecke algebras with unequal parameters in terms of equivariant homology.

For any $(\r,u)\in\Ph^K(G_\CC^*)$ we denote by $\Xi_{\r,u}$ the set of irreducible 
representations (up to isomorphism) of the group of connected components of the 
simultaneous centralizer of $(\r,u)$ in $G_\CC^*$ on which the action of the centre of 
$G_\CC^*$ is trivial.

According to \cite{\LQ}, for any $(\r,u)\in\Ph^K_1(G_\CC^*)$ the fibre of the map (b) at 
$(\r,u)$ is in bijection with $\Xi_{\r,u}$. This suggests that more generally for any
$(\r,u)\in\Ph^K(G_\CC^*)$, the fibre of the (conjectural) map (a) at $(\r,u)$ is in 
bijection with $\Xi_{\r,u}$. 

Note that in general neither side of (a) is well understood. But recent results of J.-L.Kim
\cite{\KI} give a classification of "supercuspidal representations" of $G_K$ (assuming that
$p$ is sufficiently large) which gives some hope that the left hand side of (a) can be 
understood for such $p$.

\subhead 11. Real groups\endsubhead
The Weil group of $\RR$ is by definition $W_\RR=\CC^*\T\Gal(\CC/\RR)$ with the group 
structure $(z_1,\t_1)(z_2,\t_2)=((-1)^{\e(\t_1)\e(\t_2)}z_1\t_1(z_2),\t_1\t_2)$ where \lb
$\e(\t)=0$ if $\t=1$ and $\e(\t)=1$ if $\t\ne1$.  We can identify $W_\RR$ with the subgroup
of the group of nonzero quaternions $a+bi+cj+dk$ generated by 
$\{a+bi;(a,b)\in\RR^2-\{0\}\}$ and by $j$. We
regard $W_\RR$ as a Lie group with two connected components.
Let $\Ph^\RR(G_\CC^*)$ be the set of all continuous homomorphism $W_\RR@>>>G_\CC^*$ whose 
image consists of semisimple elements, up to conjugation by $G_\CC^*$. Let $\ck$ be a 
maximal compact subgroup subgroup of the Lie group $G_\RR$. An "irreducible admissible" 
representation of $G_\RR$ is by definition a $\CC$-vector space $E$ with an action of $\ck$
and one of $\Lie(G_\RR)$ such that any vector in $E$ is contained in a finite 
dimensional $\ck$-stable subspace of $E$; the two actions induce the same action on 
$\Lie(\ck)$; the action of $\Lie(G_\RR)$ is compatible with the $\ck$-action on 
$\Lie(G_\RR)$ and the $\ck$-action on $E$. Moreover, this should be irreducible in the
obvious sense.

According to Langlands \cite{\LAA} there is a canonical finite to one surjective map
$$\{\text{irred. admissible representations of $G_\RR$ up to isom.}\}@>>>
\Ph^\RR(G_\CC^*).\tag a$$
The fibres of (a) have been described by Knapp and Zuckerman.

We now give some examples of elements of $\Ph^\RR(G_\CC^*)$ in the case where
$G_\CC^*=GL(V)$ with $V$ a finite dimensional $\CC$-vector space. Assume that we are given
a direct sum decomposition $V=\op_{p,q}V^{p,q}$ with $(p,q)\in\CC\T\CC,p-q\in\ZZ$; assume 
also that we are given a $\CC$-linear isomorphism $\ph:V@>>>V$ such that $\ph^2=1$ and 
$\ph(V^{p,q})=V^{q,p}$ for all $p,q$. We define an action of $W_\RR$ on $V$ by specifying 
the action of $(z,1)$ with $z\in\CC^*$ and the action of $(1,\t)$ with 
$\t\in\Gal(\CC/\RR)-\{1\}$. If $(p,q)$ runs only in $\ZZ\T\ZZ$ we have 
$(z,1)\cdot x=z^p\bar z^qx$ for $z\in\CC^*,x\in V^{p,q}$. The same formula holds in the 
general case: we interpret $z^p\bar z^q$ as 
$(z\bar z)^{(p+q)/2}(\fra{z}{\sqrt{z\bar z}})^{p-q}$. (The strictly positive real number 
$z\bar z$ can be raised to any complex power.) We have 
$(1,\t)\cdot x=\sqrt{(-1)^{p-q}}\ph(x)$ for $x\in V^{p,q}$. This is an object of 
$\Ph^\RR(G_\CC^*)$.

\subhead 12. Global fields\endsubhead
Let $k$ be a field which is a finite algebraic extension of the field of rational 
functions in one variable over the finite field $\FF_p$. Let $A$ be the ring of adeles of 
$k$ and let $k@>>>A$ be the canonical imbedding. Let $\bar k$ be an algebraic closure of 
$k$. Let $l$ be a prime number $\ne p$.

The global Langlands conjecture \cite{\LA} predicts a connection between the set consisting
of irreducible "cuspidal" representations of $G_A$ with nonzero vectors fixed by $G_k$ on 
the one hand and a certain set of homomorphisms of $\Gal(\bar k/k)$ into $G_{\bbq}^*$ which
are irreducible in a suitable sense, on the other hand.

This conjecture has been proved in the case where $G_k=GL_n(k)$. (See \cite{\DR}
for $n=2$ and \cite{\LF} for any $n$.)

There is an analogous conjecture in which $k$ is replaced by a finite extension of $\QQ$
and also a geometric analogue of the conjecture in which the curve over a finite field 
represented by $k$ is replaced by a smooth projective curve over $\CC$. (See \cite{\KW}.)

\subhead 13. Cells in affine Weyl groups and unipotent classes \endsubhead
Define $\ph:\ch@>>>\ca\ot J$ in terms of $\car$ as in \S4.
Let $K$ be an algebraic closure of the field $\CC(v)$ of rational functions with 
coefficients in $\CC$ in an indeterminate $v$.
Let $c$ be a two-sided cell of $\tW$. Let $J_c$ be the corresponding direct summand of the
ring $J$. We can find some simple module $E$ of the $\CC$-algebra $\CC\ot J_c$. It is 
necessarily of finite dimension over $\CC$. We can regard $K\ot_\CC E$ as a 
$K\ot_\QQ J$-module in
which the summands $K\ot_\QQ J_{c'}$ act as zero for $c'\ne c$. For $y,y'$ in $Y^+$ we 
have $T_{a^y}T_{a^{y'}}=T_{a^{y'}}T_{a^y}=T_{a^{y+y'}}$ (see \S3). Hence the operators 
$\ph(T_{a^y}):K\ot_\CC E@>>>K\ot_\CC E$
(with $y\in Y^+$) commute. We can find $e\in K\ot_\CC E-\{0\}$ which is a simultaneous 
eigenvector for these operators. Thus we have $\ph(T_{a^y})e=b(y)e$ for all $y\in Y^+$ 
where $b(y)\in K^*$ satisfy $b(y)b(y')=b(y+y')$ for any $y,y'$ in $Y^+$. There is a unique 
element $t\in K^*\ot X$ such that, if $t=\sum_sk_s\ot x_s$ with $k_s\in K^*,x_s\in X$ then 
$b(y)=\prod_sk_s^{\la y,x_\ra}$ for any $y\in Y^+$. One can show \cite{\LK} that $t$ is a 
very
special element of $K^*\ot T$: we can write uniquely $t=t't''$ where $t''\in C^*\ot X$ and
$t'\in\{v^n;n\in\ZZ\}\ot X\sub T^*_K\sub G^*_K$ is equal to 
$\ph\left(\sm v&0\\0&v\i\esm\right)$ for some homomorphism of algebraic groups 
$\ph:SL_2(K)@>>>G_K^*$. Let $C$ be the conjugacy class in $G_\CC^*$ such that
$\ph\left(\sm 1&1\\0&1\esm\right)$ is conjugate in $G_K^*$ to some element of $C$. One can
show \cite{\LK} that $C$ is well defined by $c$ (it is independent of the choice of 
$E$,$e$,$\ph$) and that $c\m C$ is a
$$\text{ bijection \{two-sided cells of }\tW\}@>\si>>
\{\text{unipotent conjugacy classes in }G^*_\CC\}.\tag a$$

\subhead 14. Special unipotent classes\endsubhead
We preserve the setup in \S13.
The intersection of $W$ with a two-sided cell of $\tW$ is said to be a two-sided cell of 
$W$ if it is nonempty. Note that the two-sided cells of $W$ form a partition of $W$. A 
unipotent conjugacy class in $G_\CC^*$ is said to be {\it special} if it corresponds under
the bijection \S13(a) to a two-sided cell of $\tW$ which has a nonempty intersection with 
$W$. 

The special unipotent classes of $G_\CC^*$ were introduced in a different (but equivalent)
way in \cite{\LB} as the unipotent classes such that the corresponding irreducible
representation of $W$ (under the Springer correspondence) is in the class $\cs_W$ defined 
in \cite{\LB}. This definition makes sense when $\CC^*$ is replaced by any algebraically
closed field $A$. For $\r\in\cs_W$ we denote by $C_{\r,A}$ the corresponding special
unipotent element of $G_A$ and by $C^*_{\r,A}$ the corresponding special unipotent element
of $G^*_A$. (The sets $\cs_W$ for $\car,\car^*$ coincide.) Let 
$\hat C_{\r,A}=\bC_{\r,A}-\cup_{C'}\bC'$ where $\bC_{\r,A}$ is the closure of $C_{\r,A}$ 
and $C'$ runs over the special unipotent classes contained in $\bC_{\r,A}-C_{\r,A}$. It is known that the subsets $\hat C_{\r,A}$ form a partition of the unipotent variety of $G_A$ 
into locally closed subvarieties which are rational homology manifolds. We define similarly
the subvarieties $\hat C^*_{\r,A}$ of the unipotent variety of $G^*_A$. We have the 
following result (see \cite{\LR},\cite{\LS}):

{\it For any $\r\in\cs_W$ there exists a polynomial $P_\r$ with integer coefficients such 
that for any $q$ (a power of a prime number) we have}
$$|\hat C_{\r,\ov{\FF}_q}\cap G_{\FF_q}|=|\hat C^*_{\r,\ov{\FF}_q}\cap G^*_{\FF_q}|
=P_\r(q).$$

\subhead 15. Preparatory results\endsubhead
Let $Y,X$ be two free abelian groups of finite rank and let $\la,\ra:Y\T X@>>>\ZZ$ be a
perfect pairing. Let $\fA:Y@>>>Y$ 
be a homomorphism such that $\det(\fA)\ne0$ that is such that
$|Y/\fA Y|<\iy$. We then have $|Y/\fA Y|=\pm\det(\fA)$. Define a homomorphism 
$\fA':X@>>>X$ by
$\la y,\fA'(x)\ra=\la\fA(y),x\ra$ for all $y\in Y,x\in X$. Then $\det(\fA')=\det(\fA)$ 
hence
$|X/\fA'(X)|<\iy$. Now $\fA$ (resp. $\fA'$) induces endomorphisms of $\QQ\ot Y$ and of 
$\QQ/\ZZ\ot Y$ (resp. $\QQ\ot X$ and $\QQ/\ZZ\ot X$) denoted again by $\fA$ (resp. $\fA'$).
Also, $\la,\ra$ induces a $\QQ$-linear 
pairing $(\QQ\ot Y)\T(\QQ\ot X)@>>>\QQ$ denoted again
by $\la,\ra$. We define a pairing $(,):Y/\fA(Y)\T X/\fA'(X)@>>>\QQ/ZZ$ by
$$(y,x)=\la\fA\i(y),x\ra\mod\ZZ=\la y,\fA'{}\i(x)\ra\mod\ZZ$$
where $y\in Y,\fA\i(y)\in\QQ\ot Y$, $x\in X,\fA'{}\i(x)\in\QQ\ot X$. Now $x\m[y\m(y,x)]$
is an isomorphism

(a) $X/\fA'(X)@>\si>>\Hom(Y/\fA(Y),\QQ/\ZZ)$.
\nl
We define an isomorphism $Y/\fA(Y)@>\si>>(\QQ/\ZZ\ot Y)^{\fA+1}$ (fixed point set of 
$\fA+1$) by 
$$y\m\text{image of $\fA\i(y)$ under }\QQ\ot Y@>>>\QQ/\ZZ\ot Y.$$
Similarly we have an isomorphism $X/\fA'(X)@>\si>>(\QQ/\ZZ\ot X)^{\fA'+1}$.
Via the last two isomorphisms, (a) becomes an isomorphism
$$(\QQ/\ZZ\ot X)^{\fA'+1}@>\si>>\Hom((\QQ/\ZZ\ot Y)^{\fA+1},\QQ/\ZZ).$$
This is induced by $\x\m[\et\m\la\fA(\et),\x\ra\mod\ZZ$ where 
$\x\in\QQ\ot X,\fA'(\x)\in X$,
$\et\in\QQ\ot Y,\fA(\et)\in Y$. Let $p$ be a prime number and let $(\QQ/\ZZ)'$ be the
subgroup of $\QQ/\ZZ$ consisting of elements of order not divisible by $p$. Assume now that
$p$ does not divide $\det(\fA)$. Then $(\QQ/\ZZ\ot Y)^{\fA+1}=((\QQ/\ZZ)'\ot Y)^{\fA+1}$ 
and we get an isomorphism
$$(\QQ/\ZZ\ot X)^{\fA'+1}@>\si>>\Hom(((\QQ/\ZZ)'\ot Y)^{\fA+1},(\QQ/\ZZ)').\tag b$$
Now let $k$ be an algebraic closure of the finite field $\FF_p$. Let $(k^*\ot Y)^{\fA+1}$ 
be the fixed point set of the endomorphism $z\ot y\m z\ot(\fA+1)y$ of $k^*\ot Y$. We define
a canonical isomorphism
$$\Hom(((\QQ/\ZZ)'\ot Y)^{\fA+1},(\QQ/\ZZ)')@>\si>>\Hom((k^*\ot Y)^{\fA+1},k^*)\tag c$$
as follows. We choose an isomorphism $\z:(\QQ/\ZZ)'@>\si>>k^*$. Then 
$\z\ot1:(\QQ/\ZZ)'\ot Y)@>\si>>k^*\ot Y$ restricts to an isomorphism
$\z_1:(\QQ/\ZZ)'\ot Y)^{\fA+1}@>\si>>(k^*\ot Y)^{\fA+1}$ and (c) carries a homomorphism
$\ph:((\QQ/\ZZ)'\ot Y)^{\fA+1}@>>>(\QQ/\ZZ)'$ to 
$\z\ph\z_1\i$. We must show that the map (c)
is independent of the choice of $\z$. Let $\k:(\QQ/\ZZ)'@>\si>>(\QQ/\ZZ)'$ be an 
isomorphism. Then $\k\ot1:(\QQ/\ZZ)'\ot Y)@>\si>>(\QQ/\ZZ)'\ot Y$ restricts to an 
isomorphism
$\k_1:(\QQ/\ZZ)'\ot Y)^{\fA+1}@>\si>>((\QQ/\ZZ)'\ot Y)^{\fA+1}$ 
and it is enough to show that
for any homomorphism $\ph:((\QQ/\ZZ)'\ot Y)^{\fA+1}@>>>(\QQ/\ZZ)'$ 
we have $\k\ph\k_1\i=\ph$.
Since $(\QQ/\ZZ)'$ is an injective $\ZZ$-module, there exists a homomorphism 
$\ti\ph:(\QQ/\ZZ)'\ot Y@>>>(\QQ/\ZZ)'$ whose restriction to $((\QQ/\ZZ)'\ot Y)^{\fA+1}$ is 
$\ph$. It is enough to show that $\k\ti\ph(\k\ot1)\i=\ti\ph$. By choosing a basis of $Y$ we
see that it is enough to show that for any homomorphism $\ps:(\QQ/\ZZ)'@>>>(\QQ/\ZZ)'$ we
have $\k\ps\k\i=\ps$. This follows from the fact that the ring of endomorphisms of the
group $(\QQ/\ZZ)'$ is commutative (it is a product of rings of $l$-adic integers for
various primes $l\ne p$).

Let $\mu$ be the group of roots of $1$ in $\CC$. We note that the isomorphism
$\z':\QQ/\ZZ@>>>\mu$ given by $r\m\exp(2\pi i r)$ induces an isomorphism 
$\z'\ot1:\QQ/\ZZ\ot X@>>>\mu\ot X$ and this restricts to an isomorphism 
$\z'_1:(\QQ/\ZZ\ot X)^{\fA'+1}@>\si>>(\mu\ot X)^{\fA'+1}$ where $(\mu\ot X)^{\fA'+1}$ is 
the 
fixed point set of the endomorphism $z\ot x\m z\ot(\fA'+1)x$
of $\mu\ot X$. Via $\z'_1$ and 
(c), the isomorphism (b) becomes a canonical isomorphism
$$(\mu\ot X)^{\fA'+1}@>\si>>\Hom((k^*\ot Y)^{\fA+1},k^*).\tag d$$

\subhead 16. Groups over $\FF_q$\endsubhead
We return to the setup in \S8(a). Let 
$k$ be an algebraic closure of the finite field $\FF_p$. Let $K$ be an algebraically closed
field of characteristic $0$ with a fixed imbedding of groups $\io:k^*@>>>K^*$. 

We have $T_k=k^*\ot Y$, $T^*_\CC=\CC^*\ot X$. Let $\FF_q$ be the subfield of $k$ such that
$|\FF_q|=q$. The ring 
homomorphism $k@>>>k,c\m c^q$ induces (as in \S6) a group homomorphism $F:G_k@>>>G_k$ 
(Frobenius map) whose fixed point set is the finite group $G_{\FF_q}$. Following
\cite{\DL} we consider for any $w\in W$ the set $\dX_w=\{g\in G_k; g\i F(g)\in\dw U_k^+\}$,
an algebraic variety over $k$. Let $T_k^w=\{t\in T_k;t^q=w\i(t)\}$, a finite subgroup of
$T_k$. The finite group $G_{\FF_q}\T T_k^w$ acts on $\dX_w$ by $(g_1,t):g\m g_1gt\i$. Let 
$\c_w:G_{\FF_q}\T T_k^w@>>>\ZZ$ be the class function which to any  
$(g_1,t)\in G_{\FF_q}\T T_k^w$ associates the alternating sum of traces of $(g_1,t)^*$ on 
the $l$-adic cohomology with compact support of $\dX_w$. (Here $l$ is any prime number 
$\ne p$ but the resulting class function is known to be independent of $l$; see 
\cite{\DL}.) For any irreducible $G_{\FF_q}$-module $\r$ over $K$ let $\ce_\r$ be the 
set of all pairs $(w,\th)$ where $w\in W$, $\th\in\Hom(T_k^w,K^*)$ and 
$\sum_{(g_1,t)\in G_{\FF_q}\T T_k^w}\th(t)\tr(g_1,\r)\c_w(g_1,t)\ne0$. According to 
\cite{\DL} we have $\ce_\r\ne\em$ for any $\r$.  

To any $(w,\th)\in\ce_\r$ we associate an element 
$\hat\th\in T_\CC^{*w\i}:=\{t\in T^*_\CC;t^q=w(t)\}$ as 
follows. Define $\fA:Y@>>>Y$ by $y\m qw(y)-y$ and $\fA':X@>>>X$ by $x\m qw\i(x)-x$. Then 
$T_k^w=(k^*\ot Y)^{\fA+1}$ is a finite group of order prime to $p$. Hence 
$\th:T_k^w@>>>K^*$ has values in the group of roots of $1$ of order prime to $p$ in $K^*$ 
which can be identified with $k^*$ via $\io$. Thus $\th$ can be viewed as an element of 
$\Hom((k^*\ot Y)^{\fA+1},k^*)$ 
so that it corresponds under \S15(d) to an element $\hat\th$ 
of $(\mu\ot X)^{\fA'+1}$. 
This last group is a subgroup of $(\CC^*\ot X)^{\fA'+1}$ (the fixed 
point set of the endomorphism $z\ot x\m z\ot(\fA'+1)x$ 
of $\CC^*\ot X$) which is the same as
$T_\CC^{*w\i}$. From the results in \cite{\DL} we see that the $W$-orbit of $\hat\th$ in 
$T_\CC^*$ depends only on $\r$ and not on the choice of $(w,\th)$ in $\ce_\r$. We thus have
a well defined map
$$\align&\{\text{irred. $G_{\FF_q}$-modules over $K$ up to isom.}\}@>>>\\&
\{\text{semisimple conjugacy classes in $G_\CC^*$ stable under $g\m g^q$}\};\tag a\endalign
$$
it is given by $\r\m G_\CC^*\text{-conjugacy class of }\hat\th$ (as above).
This map appears in \cite{\DL} in a somewhat different form. In \cite{\DL} $G_\CC^*$ is
replaced by $G^*_k$. But the method of \cite{\DL} is less canonical: it is based on two 
choices (see \cite{\DL, (5.0.1), (5,0.2)}) while the present method is based on only one 
choice, that of $\io$; the choice of $\io$ can be also eliminated as we will see below).

An element $g\in G_\CC^*$ is said to be special if the unipotent part $g_u$ of $g$ is a 
special unipotent element (see \S14)
of the connected centralizer of the semisimple part $g_s$ of 
$g$ (a reductive connected group). A conjugacy class in $G_\CC^*$ is said to be special if 
some/any element of it is special. The map (a) can be refined to a canonical map
$$\align&\{\text{irred. $G_{\FF_q}$-modules over $K$ up to isom.}\}@>>>\\&
\{\text{special conjugacy classes in $G_\CC^*$ stable under }g\m g^q\}.\tag b\endalign$$
(See \cite{\LI}.) Note that (a) is the composition of (b) with the map which to the 
$G_\CC^*$-conjugacy class of a special element $g$ associates the $G_\CC^*$-conjugacy class
of $g_s$. The map (b) is surjective and its fibres are described explicitly in \cite{\LI},
\cite{\LM}.

Note that the maps (a),(b) depend on the choice of the imbedding $\io:k^*@>>>K^*$. 
However if we take $K$ to be an algebraic closure of the quotient field of the ring of Witt
vectors of $k$ then there is a canonical choice of $\io$ and the maps (a),(b) become
completely canonical.

\subhead 17. Character sheaves\endsubhead
Define $B_\CC^+,U_\CC^+,\dw$ in terms of $\CC,\car$ as in \S6. Let $\ce$ be a $\CC$-local 
system of rank $1$ on $T_\CC$ with finite monodromy. The monodromy of $\ce$ is a 
homomorphism $f:Y@>>>\CC^*$ with finite image which can be viewed as an element of finite 
order $\c_\ce\in\CC^*\ot X$ given by $\c_\ce=\sum_jf(y_j)\ot x_j$ where $(y_j)$ is a basis
of $Y$ and $(x_j)$ is the dual basis of $X$. Moreover $\ce\m\c_\ce$ is a bijection
$$\align& \{\CC-\text{local systems of rank $1$ on $T_\CC$ with finite monodromy 
up to isom.}\}@>\si>>\\&
\{\text{elements of finite order of }T_\CC^*\}.\tag a\endalign$$
Let $c:G_\CC@>>>G_\CC/U_\CC^+$ be the obvious map. An irreducible intersection cohomology 
complex $K$ on $G_\CC$ is said to be a {\it character sheaf} on $G_\CC$ if it is 
$G_\CC$-equivariant and if for any $w\in W$ and any $j\in\ZZ$ the $j$-th cohomology sheaf 
of $c_!K$ restricted to $B_\CC^+\dw B_\CC^+/U_\CC^+$ is a local system $\cl_{K,w,j}$ with 
finite monodromy. 
We can find $w,j$ as above and a local system $\ce$ of rank $1$ on $T_\CC$ with finite
monodromy such that $\ce$ is a direct summand of the inverse image of $\cl_{K,w,j}$ under
the map $T_\CC@>>>B_\CC^+\dw B_\CC^+/U_\CC^+$, $t\m\dw tU_\CC^+$. 
One can show that the corresponding
element $\c_{\ce}\in T^*_\CC$ is well defined (up to the action of $W$) that is, it
does not depend on the choice of $w,j,\ce$. Thus we have a well defined map
$$\align&\{\text{character sheaves on $G_\CC$ up to isom.}\}@>>>\\&
\{\text{conjugacy classes of elements of finite order in }G_\CC^*\};\tag a\endalign$$
it is given by $K\m G_\CC^*\text{-conjugacy class of }\c_{\ce}$ (as above). The 
map (a) can be refined to a canonical map
$$\align&\{\text{character sheaves on $G_\CC$ up to isom.}\}@>>>\\&
\{\text{special conjugacy classes in $G_\CC^*$ of elements $g$ with $g_s$ of finite order}
\}.\tag b\endalign$$
(See \cite{\LJ}.) Note that (a) is the composition of (b) with the map which to the 
$G_\CC^*$-conjugacy class of a special element $g$ associates the $G_\CC^*$-conjugacy class
of $g_s$. The map (b) is surjective and its fibres are described explicitly in \cite{\LJ}.

\subhead 18. Vogan duality\endsubhead 
To $G_\CC$ we associate a finite collection of polynomials: 
those recording the restrictions
of the cohomology sheaves of the simple perverse sheaves on $\cb_\CC$, equivariant under 
the conjugation action of the centralizers of the various involutions of $G_\CC$, 
to the various orbits of these centralizers. (These polynomials were studied in 
\cite{\LV}.) We consider also the analogous collection of polynomials associated to 
$G^*_\CC$. {\it Vogan duality} \cite{\VO} states that these two collections of polynomials
are related to each other by a simple algebraic rule: essentially the inversion of a 
matrix. This is a generalization of the inversion formula in \cite{\KL}.

\subhead 19. Multiplicities in standard modules\endsubhead 
Let $K$ be either as in \S10 or $K=\RR$. The Grothendieck group whose basis consists of 
admissible irreducible representations of $G_K$ has another
natural basis consisting of "standard representations" in natural bijection with the first basis.
The representations in the second 
basis are easier to describe and understand. Hence the (upper triangular) matrix $\cm$ 
expressing the second basis in terms of the first basis ("multiplicity matrix") is of interest. In 
every known case the entries of $\cm$ can be expressed in terms of intersection 
cohomology coming from the geometry of $G_\CC^*$. For the case where $K=\RR$ we refer the
reader to \cite{\ABV}; in this case some of the polynomials in \S18 (attached to
$G_\CC^*$) evaluated at $1$ appear as entries of $\cm$. In the remainder of this 
subsection we assume that $K$ is as in \S10. For simplicity we assume that $\car$ is of adjoint
type. We use the notation in \S10. 

We fix an element $w^0\in W_K$ such that $\o(w^0)=1$. 
Let $\ti\Ph^K$ be the set of all triples $(\r,u,E)$ (up to $G_\CC^*$-conjugacy) where 
$\r:W_K@>>>G_\CC^*$ is an admissible homomorphism, $u$ is a unipotent element of $G^*_\CC$
such that $\r(w)u\r(w)\i=u^{q^{\o(w)}}$ for all $w\in W_K$
and $E$ is an irreducible representation of the group of connected components of
$G_{\r,u}^*:=\{g\in G_\CC^*;\r(w)g\r(w)\i=g\text{ for all }w\in W_K,gu=ug\}$ on which the 
image of the centre of $G_\CC^*$ acts trivially.

By \S10 it is expected that $\ti\Ph^K$ is an index set for both the rows and the columns
of $\cm$. We shall describe a matrix $\cm'$ indexed by $\ti\Ph^K$ which is defined in 
terms of geometry of $G_\CC^*$.

Let $\Ps$ be the set of homomorphisms $\ps:\ci@>>>G^*_\CC$ such that $\ps(\ci)$ is finite
and such that $\G_\ps:=\{g\in G^*_\CC;g\ps(w)g\i=\ps(w^0ww^0{}\i)\text{ for all }w\in\ci\}$
is non-empty. Let $\bar\Ps$ be the set of $G^*_\CC$-orbits (by conjugacy) on $\Ps$.

Define $\k:\ti\Ph^K@>>>\bar\Ps$ by $(\r,u,E)\m\r|_\ci$. 
The entries $m_{(\r,u,E),(\r',u',E')}$ of $\cm'$ can be described as follows. If 
$(\r,u,E),(\r',u',E')$ are not in the same fibre of $\k$ then $m_{(\r,u,E),(\r',u',E')}=0$.

We now fix $\ps\in\Ps$.
Let $G^*_\ps=\{g\in G^*_\CC;g\ps(w)g\i=\ps(w)\text{ for all }w\in\ci\}$. This is the 
centralizer of a finite subgroup of $G_\CC^*$ hence is a (possibly disconnected) reductive
subgroup of $G^*_\CC$. Let $G'_\ps{}^*$ be the normalizer of $\ps(\ci)$ in $G^*_\CC$. Note
that $G^*_\ps$ is a normal subgroup of finite index of $G'_\ps{}^*$ and that $\G_\ps$ is a
single $G^*_\ps$-coset in $G'_\ps{}^*$. The
fibre of $\k$ at $\ps$ can be identified with the set of all triples $(s,N,E)$ (up to 
$G^*_\ps$-conjugacy) where $s\in G'_\ps{}^*$ is a semisimple element such that 
$s\in\G_\ps$, $N$ is an element of $X_s:=\{N_1\in\Lie(G^*_\ps);\Ad(s)N_1=qN\}$ (necessarily 
nilpotent) and $E$ is an irreducible representation of the group of connected components of
$G^*_{\ps,s,N}:=\{g\in G^*_\ps;gs=sg,\Ad(g)N=N\}$ on which the image of the centre of 
$G_\CC^*$ acts trivially. (The identification is given by $(\r,u,E)\m(\r(w^0),\log(u),E)$ 
where $\r|_\ci=\ps$.) We now consider two elements $(s,N,E),(s',N',E')$ in $\k\i(\ps)$. If
$s,s'$ are not in the same $G^*_\ps$-orbit then $m_{(s,N,E),(s',N',E')}=0$. Now assume that
$(s,N,E),(s',N',E')$ are such that $s,s'$ are in the same $G^*_\ps$-orbit. We can assume 
that $s=s'$. Let $G^*_{\ps,s}=\{g\in G^*_\ps;gs=sg\}$. This is an algebraic group which 
acts on $X_s$ by conjugation with finitely many orbits. Let $C$ be the $G^*_{\ps,s}$-orbit
of $N$ and let $C'$ be the $G^*_{\ps,s}$-orbit of $N'$.
Note that $E$ (resp. $E'$) determines a local system $\uE$ (resp. $\uE'$) on $C$ (resp. 
$C'$) which is $G^*_{\ps,s}$-equivariant and is irreducible as a $G^*_{\ps,s}$-equivariant
local system. If $C$ is not contained in the closure of $C'$ then 
$m_{(s,N,E),(s',N',E')}=0$. Now assume that $C$ is contained in the closure of $C'$. Let 
$\uE'{}^\sh$ be the intersection cohomology complex on the closure of $C'$ determined by
$\uE'$. For every integer $j$ we consider the $j$-th cohomology sheaf of $\uE'{}^\sh$ 
restricted to $C$; this is a $G^*_{\ps,s}$-equivariant local system on $C$ in which $\uE$ 
appears say $n_j$ times. We set $m_{(s,N,E),(s',N',E')}=\sum_j(-1)^jn_j$.

We see that the intersection cohomology complexes on $X_s$ considered above are essentially
of the type considered in \cite{\LU}. We also see that the objects in $\k\i(\ps)$ behave 
like the parameters for the unipotent representations for a collection of not necessarily 
split and not necessarily connected $p$-adic groups smaller than $G_K$.

It is known that $\cm$, $\cm'$ coincide as far as the entries 
with both indices contained in $\k\i(1)$ are concerned; these corresponds to unipotent representations.
(This was conjectured by the author and independently, in a special case connected with $GL_n$, in \cite{\ZE}; 
the proof was given by Ginzburg in a special case connected with the affine Hecke algebra $\ch$ 
and by the author in the general case.) We expect that $\cm=\cm'$. (See also \cite{\VOO}.)

\subhead 20. Multiplicities in tensor products\endsubhead
Assume that $\car$ is simply connected. For $\l,\l',\l''$ in $X^+$ let $m_{\l,\l',\l''}$ be
the multiplicity of $\L_{\l''}$ in the tensor product $\L_\l\ot\L_{\l'}$ (an object of 
$\cc$, see \S5). On the other hand let $l^*:\tW^*@>>>\NN,T_w,c_w,p_{w,z}$ be defined
like $l:\tW@>>>\NN,T_w,c_w,p_{w,z}$ in \S3,\S4 but with respect to $\car^*$ instead of 
$\car$. We have $\tW^*=\{wa^x;w\in W,x\in X\}$. 

For any $\l\in X^+$ there is a unique
element $M_\l$ in the double coset $Wa^\l W$ on which $l^*:Wa^\l W@>>>\NN$ achieves its 
maximum value. For any $\l,\l'$ in $X^+$ we have in $\ch^*$:
$$(P\i c_{M_\l})(P\i c_{M_{\l'}})=\sum_{\l''\in X^+}\tm_{\l,\l',\l''}(P\i c_{M_{\l''}})$$
where $P\in\ca$ is given by $c_{M_0}c_{M_0}=Pc_{M_0}$ and $\tm_{\l,\l',\l''}\in\ca$. In 
\cite{\LE} it is shown that
$$m_{\l,\l',\l''}=\tm_{\l,\l',\l''};\tag a$$
in particular, 
$$\tm_{\l,\l',\l''}\text{ is a constant}.\tag b$$
(In the special case where $G_\CC$ is a general linear group, this was proved earlier in 
\cite{\LD} using the theory of Hall-Littlewood functions.) In \cite{\LE} it is also shown 
that for $\l,\l''$ in $X^+$, 

(c) {\it $\dim\L_\l^{\l'}$ (that is the multiplicity of the weight $\l'$ in the 
$\uU$-module $\L_\l$) is equal to $p_{M_{l'},M_\l}(1)$.}
\nl
Note that at the time when \cite{\LE} was written it was known that the product of elements
of the form $P\i c_{M_\l}$ corresponds to the convolution of $G^*_{\CC[[\e]]}$-equivariant
simple perverse sheaves on the "affine Grassmannian" $G^*_{\CC((\e))}/G^*_{\CC[[\e]]}$ so 
that (b) is equivalent to the statement that such a
convolution is a direct sum of simple perverse sheaves of the same type (without shift).
Thus it was clear that the category whose objects are finite direct sums of 
$G^*_{\CC[[\e]]}$-equivariant simple perverse sheaves on the "affine Grassmannian" has a
natural monoidal structure given by convolution; moreover (b) showed that this monoidal
category was very similar to that of representations of $G_\CC$ (identical at the level of
Grothendieck groups). But it was not clear how to construct the commutativity isomorphism
for the convolution product. This was accomplished around 1989 by V.Ginzburg \cite{\GI} and
later in a more elegant form by V.Drinfeld. As a result, $G_\CC$ can be reconstructed from
the tensor category of $G^*_{\CC[[\e]]}$-equivariant perverse sheaves on 
$G^*_{\CC((\e))}/G^*_{\CC[[\e]]}$, see \cite{\GI}.

\subhead 21. Canonical bases\endsubhead
Define $\uf$ as in \S5 in terms of $\car$. Let $A=\RR[[\e]]$ where $\e$ is an 
indeterminate. Define $U^{*+}_A$ in terms of $\car^*$ in the same way as $U^+_A$ was 
defined in \S6 in terms of $\car$. By \cite{\LP,\S10} there is a canonical bijection 
between the {\it canonical basis} of $\uf$ (defined as in \cite{\LN}, \cite{\LO}) and a 
certain collection of subsets of $U^{*+}_A$ which form a partition of the totally positive
part of $U^{*+}_A$. The bijection is not defined directly; instead it is shown that both 
sets are parametrized by the same combinatorial objects.

\subhead 22. Modular representations\endsubhead
Let $k$ be an algebraic closure of the finite field $\FF_p$. Assume that $\car$ is simply 
connected. For $\l\in X^+$ the $G_k$-module $\L_{\l,k}$ (see \S7) is not necessarily
irreducible but has a unique irreducible quotient $\L_{\l,k}^\sh$. For $\l,\l'$ in $X^+$ 
let $m_{l,\l'}$ be the number of times that $\L_{\l',k}^\sh$ appears in a composition 
series of the $G_k$-module $\L_{\l,k}$. Note that the knowledge of the multiplicities 
$m_{\l,\l'}$ implies the knowledge of the character of the $G_k$-modules $\L_{\l,k}^\sh$ 
since the character of $\L_{\l,k}$ is known by Weyl's character formula. Conjecturally (see
\cite{\LC}) if $p$ is sufficiently large with respect to $\car$, the multiplicities 
$m_{\l,\l'}$ can be expressed in terms of polynomials $p_{w,z}$ (as in \S4) where $w,z$ are
elements in $\tW^*$ which have maximal length in their left $W$-coset; they can be also 
expressed in terms of certain intersection cohomology spaces associated with the geometry 
of $G^*_{\CC((\e))}$, where $\e$ is an indeterminate. A proof of the conjecture (without an
explicit bound for $p$) is provided by combining \cite{\AJS}, \cite{\KT}, \cite{\KLLL} or 
alternatively by combining \cite{\AJS}, \cite{\ABG}.


\widestnumber\key{ABG}
\Refs
\ref\key{\ABV}\by J.Adams, D.Barbasch and D.A.Vogan, Jr.\paper The Langlands classification of
irreducible characters of real reductive groups, Progress in Math\publ Birkhauser\yr1992
\endref
\ref\key{\AJS}\by H.Andersen, W.Soergel and J.Jantzen\paper Representations of quantum 
groups at a $p$-th root of unity and of semisimple groups in characteristic $p$: 
independence of $p$\jour Ast\'erisque\vol220\yr1994\endref
\ref\key{\ABG}\by S.Arkhipov, R.Bezrukavnikov and V.Ginzburg\paper Quantum groups, the loop
grassmannian and the Springer resolution\jour Jour.Amer.Math.Soc.\vol17\yr2004\pages595-678
\endref
\ref\key{\CH}\by C.Chevalley\paper Certains sch\'emas de groupes semi-simples\inbook
S\'em. Bourbaki 1960/61\publ Soc. Math. France\yr1995\endref
\ref\key{\DL}\by P.Deligne and G.Lusztig\paper Representations of reductive groups over
finite fields\jour Ann. Math.\vol103\yr1976\pages103-161\endref
\ref\key{\DR}\by V.Drinfeld\paper Two dimensional $l$-adic representations of the
fundamental group of a curve over a finite field and automorphic forms on $GL(2)$\jour
Amer.Jour.Math.\vol105\yr1983\pages85-114\endref
\ref\key{\GI}\by V.Ginzburg\paper Perverse sheaves on a loop group and Langlands duality
\jour math.AG/9511007\endref
\ref\key{\HT}\by M.Harris and R.Taylor\book On the geometry and cohomology of some simple
Shimura varieties\bookinfo Ann.Math.Studies\vol151\publ Princeton Univ.Press\yr2001\endref
\ref\key{\IM}\by N.Iwahori and H.Matsumoto\paper On some Bruhat decomposition and the
structure of the Hecke ring of $p$-adic Chevalley groups\jour Publ.Math.IHES\vol25\yr1965
\pages5-48\endref
\ref\key{\KW}\by A.Kapustin and E.Witten\paper Electric-magnetic duality and the geometric
Langlands program\jour arXiv:hep-th/0604151\endref
\ref\key{\KT}\by M.Kashiwara and T.Tanisaki\paper The Kazhdan-Lusztig conjecture for affine
Lie algebras with negative level\jour Duke Math.J.\vol77\yr1995\pages21-62\endref
\ref\key{\KL}\by D.Kazhdan and G.Lusztig\paper Representations of Coxeter groups and Hecke
algebras\jour Invent.Math.\vol53\yr1979\pages165-184\endref
\ref\key{\KLL}\by D.Kazhdan and G.Lusztig\paper Proof of the Deligne-Langlands conjecture 
for Hecke algebras\jour Invent.Math.\vol87\yr1987\pages153-215\endref
\ref\key{\KLLL}\by D.Kazhdan and G.Lusztig\paper Tensor structures arising from affine Lie 
algebras, IV\jour Jour. Amer. Math. Soc.\vol7\yr1994\pages383-453\endref
\ref\key{\KI}\by J.-L.Kim\paper Supercuspidal representations: an exhaustion theorem\jour
Jour.Amer.Math.Soc.\vol20\yr2007\pages273-320\endref
\ref\key{\KO}\by B.Kostant\paper Groups over $\ZZ$\inbook Algebraic Groups and Their 
Discontinuous Subgroups\bookinfo Proc. Symp. Pure Math.\vol8 publ. Amer. Math.
Soc.\yr1966
\pages90-98\endref
\ref\key{\LF}\by L.Lafforgue\paper Chtoukas de Drinfeld et correspondance de Langlands\jour
Invent.Math.\vol147\yr2002\pages1-242\endref
\ref\key{\LA}\by R.P.Langlands\paper Problems in the theory of automorphic forms\inbook
Lectures in Modern Analysis and Applications\bookinfo Lecture Notes in Math\vol170\publ
Springer Verlag\yr1970\pages18-61\endref
\ref\key{\LAA}\by R.P.Langlands\paper On the classification of irreducible representations of
real algebraic groups\inbook Representation theory and harmonic analysis on semisimple Lie groups
\pages101-170\bookinfo Math. Surveys Monogr.\vol31 Amer. Math. Soc.\publaddr Providence, RI\yr1989
\endref
\ref\key{\LB}\by G.Lusztig\paper A class of irreducible representations of a Weyl group
\jour Proc. Kon. Nederl. Akad. (A)\vol82\yr1979\pages323-335\endref
\ref\key{\LC}\by G.Lusztig\paper Some problems in the representation theory of finite 
Chevalley groups\jour Proc. Symp. Pure Math. Amer. Math. Soc.\vol37\yr1980\pages313-317
\endref
\ref\key{\LD}\by G.Lusztig\paper Green polynomials and singularities of unipotent classes
\jour Adv. Math.\vol42\yr1981\pages169-178\endref
\ref\key{\LE}\by G.Lusztig\paper Singularities, character formulas and a $q$-analog of 
weight multiplicities\jour Ast\'erisque\vol101-102\yr1983\pages208-229\endref
\ref\key{\LG}\by G.Lusztig\paper Some examples of square integrable representations of 
semisimple p-adic groups\jour Trans.Amer.Math.Soc.\vol227\yr1983\pages623-653\endref
\ref\key{\LI}\by G.Lusztig\book Characters of reductive groups over a finite field
\bookinfo Ann.Math.Studies\vol107\publ Princeton Univ.Press\yr1984\endref
\ref\key{\LII}\by G.Lusztig\paper Equivariant $K$-theory and representations of Hecke 
algebras\jour Proc. Amer. Math. Soc.\vol94\yr1985\pages337-342\endref
\ref\key{\LJ}\by G.Lusztig\paper Character sheaves, V\jour Adv.Math.\vol 61\yr1986\pages
103-155\endref
\ref\key{\LK}\by G.Lusztig\paper Cells in affine Weyl groups, II\jour J.Alg.\vol109\yr1987
\pages536-548\moreref IV\jour J. Fac. Sci. Tokyo U.(IA)\vol36\yr1989\pages297-328\endref
\ref\key{\LM}\by G.Lusztig\paper On representations of reductive groups with disconnected 
center\jour Ast\'erisque\vol168\yr1988\pages157-166\endref
\ref\key{\LN}\by G.Lusztig\paper Canonical bases arising from quantized enveloping 
algebras\jour Jour. Amer. Math. Soc.\vol3\yr1990\pages447-498\endref
\ref\key{\LNN}\by G.Lusztig\paper Affine quivers and canonical bases\jour Publ.Math.IHES
\vol76\yr1992\pages111-163\endref
\ref\key{\LO}\by G.Lusztig\book Introduction to quantum groups\bookinfo Progress in Math.
\vol110\publ Birkhauser\yr1993\endref
\ref\key{\LP}\by G.Lusztig\paper Total positivity in reductive groups\inbook Lie theory and
geometry, Progr.in Math.\vol123\publ Birkh\"auser Boston\yr1994\pages531-568
\endref
\ref\key{\LQ}\by G.Lusztig\paper Classification of unipotent representations of simple 
$p$-adic groups\jour Int. Math. Res. Notices\yr1995\pages517-589\moreref II\jour Represent.
Th.\vol6\yr2002\pages243-289\endref
\ref\key{\LR}\by G.Lusztig\paper Notes on unipotent classes\jour Asian J.Math.\vol1\yr1997
\pages194-207\endref
\ref\key{\LS}\by G.Lusztig\paper Unipotent elements in small characteristic\jour Transform.
Groups\vol10\yr2005\pages449-487\moreref II,arXiv:RT/0612320\endref
\ref\key{\LU}\by G.Lusztig\paper Graded Lie algebras and intersection cohomology,
arXiv:RT/0604535\endref
\ref\key{\LT}\by G.Lusztig\paper Study of a $\ZZ$-form of the coordinate ring of a 
reductive group\jour arxiv:0709.1286\endref
\ref\key{\LV}\by G.Lusztig and D.A.Vogan, Jr.\paper Singularities of closures of $K$-orbits on a 
flag manifold\jour Invent.Math.\vol71\yr1983\pages365-379\endref
\ref\key{\MK}\by J.McKay\paper Graphs, singularities and finite groups
\jour Proc. Symp. Pure Math. Amer. Math. Soc.\vol37\yr1980\pages183-186\endref
\ref\key{\SL}\by P.Slodowy\book Simple singularities and simple algebraic groups\bookinfo 
Lecture Notes in Math.\vol815\publ Springer Verlag\yr1980\endref
\ref\key{\VO}\by D.A.Vogan, Jr.\paper Irreducible characters of semisimple Lie groups, IV:
character multiplicity duality\jour Duke Math.J.\vol4\yr1982\pages943-1073\endref
\ref\key{\VOO}\by D.A.Vogan, Jr.\paper The local Langlands conjecture\inbook Representation
theory of groups and algebras\pages305-379\bookinfo Contemp. Math.\vol145 Amer. Math. Soc.
\publaddr Providence, RI\yr1993\endref
\ref\key{\XI}\by N.Xi\paper Representations of affine Hecke algebras and based rings of
affine Weyl groups\jour Jour.Amer.Math.Soc.\vol20\yr2007\pages211-217\endref
\ref\key{\ZE}\by A.Zelevinsky\paper A $p$-adic analogue of the Kazhdan-Lusztig conjecture
\jour Funkt.Anal.Pril.\vol15\yr1981\pages9-21\endref
\endRefs
\enddocument